\documentclass[11pt]{amsart}

\usepackage[a4paper,margin=1.05in]{geometry}
\usepackage{amsmath,amssymb,amsthm,mathtools}
\usepackage[expansion=false]{microtype}
\usepackage{booktabs,tabularx,array}
\usepackage[colorlinks=true,linkcolor=blue!50!black,citecolor=blue!50!black,urlcolor=blue!50!black]{hyperref}
\usepackage{xcolor}
\usepackage[nameinlink,noabbrev]{cleveref}

\numberwithin{equation}{section}
\allowdisplaybreaks

\crefname{lemma}{lemma}{lemmas}
\Crefname{lemma}{Lemma}{Lemmas}
\crefname{proposition}{proposition}{propositions}
\Crefname{proposition}{Proposition}{Propositions}
\crefname{definition}{definition}{definitions}
\Crefname{definition}{Definition}{Definitions}
\hypersetup{pdftitle={Convex entire solutions of the k-Hessian equation}}

\numberwithin{equation}{section}
\newtheorem{theorem}{Theorem}[section]
\newtheorem{proposition}[theorem]{Proposition}
\newtheorem{lemma}[theorem]{Lemma}

\theoremstyle{definition}
\newtheorem{definition}[theorem]{Definition}
\theoremstyle{remark}

\newcommand{\R}{\mathbb R}
\newcommand{\Sym}{\operatorname{Sym}}
\newcommand{\tr}{\operatorname{tr}}
\newcommand{\diag}{\operatorname{diag}}
\newcommand{\dist}{\operatorname{dist}}
\newcommand{\Id}{\mathrm I}
\newcommand{\PP}{\mathcal P}
\newcommand{\MA}{\operatorname{MA}}
\newcommand{\norm}[1]{\left\lVert#1\right\rVert}
\newcommand{\abs}[1]{\left|#1\right|}
\newcommand{\dd}{\,dx}

\title{Rigidity for Convex entire solutions to the \(k\)-Hessian equation}

\author{Xingchen Zhou \& Ruixuan Zhu}

\address{Xingchen Zhou, 
School of Mathematics and Statistics, Hainan University, Haikou, 570228, PR China.}
\email{zxc3zxc4zxc5@stu.xjtu.edu.cn}

\address{Ruixuan Zhu, 
Institute for Theoretical Sciences, Westlake University, Hangzhou, 310030, China.}
\email{zhuruixuan@westlake.edu.cn}

\date{}

\begin{document}

\begin{abstract}
We prove that for each \(3\le k\le n-1\) a convex \(C^2\) entire solution of
\(\sigma_k\) equation in \(\R^n\) with positive constant right hand side  is a quadratic polynomial,
 without a growth assumption on the solution itself.

\end{abstract}


\maketitle

\section{Introduction}

Let \(\sigma_k\) denote the \(k\)-th elementary symmetric
polynomial of the eigenvalues of a symmetric matrix.
We prove that every convex entire solution of
\begin{equation}\label{eq:main}
 \sigma_k(D^2u)=c_0>0\qquad\text{in }\R^n
\end{equation}
is a quadratic polynomial.

\begin{theorem}\label{thm:main}
Let \(n\ge4\), \(3\le k\le n-1\), and \(c_0>0\).
If \(u\in C^2(\R^n)\) satisfies \(D^2u\ge0\) and
\cref{eq:main}, then
\[
 u(x)=\frac12 x^TQx+a\cdot x+b
\]
for a constant matrix \(Q\ge0\), \(a\in\R^n\), and \(b\in\R\).
\end{theorem}

For \(k=1\), the equation is
\(\Delta u=c_0\). Convexity bounds every entry of \(D^2u\), and
the Liouville theorem for harmonic functions implies that the
Hessian is constant. For \(k=n\), the same conclusion is the
J\"orgens--Calabi--Pogorelov theorem for the Monge--Amp\`ere
equation \cite{Jorgens,Calabi,Pogorelov}; see also the proof of
Cheng and Yau \cite{ChengYau}.
The case \(k=2\) follows from earlier results by Shankar and Yuan \cite{SY} under semiconvex assumptions.

\vspace{0.7\baselineskip}

Recall the elliptic cone
\[
 \Gamma_k=\{\lambda\in\R^n:\sigma_j(\lambda)>0
                 \text{ for }1\le j\le k\}.
\]
A function is \(k\)-admissible, or \(k\)-convex, when
\(\lambda(D^2u)\in\Gamma_k\); see
Caffarelli, Nirenberg, and Spruck \cite{CNS}.
Every convex solution of \cref{eq:main} is \(k\)-admissible.
We also call a function semiconvex if \(D^2u\ge-K\Id\) for
some constant \(K\ge0\). 

\vspace{0.7\baselineskip}

Bao, Chen, Guan, and Ji \cite{BCGJ} proved rigidity for strictly
convex entire solutions under a quadratic lower bound at infinity:
for some \(a>0\) and \(b\in\R\),
\begin{equation}\label{eq:quadraticgrowth}
 u(x)\ge a|x|^2-b\qquad\text{for all sufficiently large }|x|.
\end{equation}
They also used the Legendre transform to pass from Hessian equations
to positive Hessian quotient equations. For the transformed solution,
the corresponding growth assumption is a quadratic upper bound.
Li, Ren, and Wang \cite{LRW} proved rigidity under
\cref{eq:quadraticgrowth} for \((k+1)\)-convex solutions.
Their proof uses a Pogorelov estimate and extends the interior
estimate approach of Chou and Wang \cite{CW}.
For locally strictly convex solutions, Du \cite{Du} gave conditions
on sublevel sets that characterize when a solution is quadratic.
These include
reverse isoperimetric, volume growth, and integral conditions.

The quadratic Hessian equation admits stronger conclusions.
Chang and Yuan \cite{CY} proved that every entire solution of
\(\sigma_2(D^2u)=1\) is quadratic if, for some \(\delta>0\),
\[
 D^2u\ge\left(\delta-\sqrt{\frac{2}{n(n-1)}}\right)\Id.
\]
This includes all convex solutions and requires no growth
assumption. Shankar and Yuan \cite{SY} subsequently proved
rigidity for every semiconvex entire solution. Their proof combines
a Jacobi inequality, the Legendre--Lewy transform, and an
approximation argument.

For general \(k\), several results allow weaker convexity while
retaining conditions at infinity. Wang and Bao \cite{WB} proved
rigidity for \(k\)-admissible solutions with two-sided quadratic
growth. The results of Chu and Dinew \cite{CD} apply to
\(k\)-admissible solutions satisfying both
\(\sigma_{k+1}(D^2u)\ge-A\) and \cref{eq:quadraticgrowth}.
Zhang \cite{Zhang} proved a concavity inequality and Pogorelov
estimates for semiconvex solutions, and obtained rigidity under
\cref{eq:quadraticgrowth} for \(2<k<n\).
Under the same lower growth bound, Tu \cite{Tu} treated all
\((n-1)\)-admissible solutions, and B.~Wang \cite{BW} treated
\(k\)-admissible solutions in the range \(2k>n\).
The principal results are compared in \cref{tab:liouville}.

\begin{table}[t]
\centering
\caption{Selected rigidity results for \(\sigma_k(D^2u)=1\) in
\(\R^n\). }
\label{tab:liouville}
\small
\renewcommand{\arraystretch}{1.18}
\setlength{\tabcolsep}{4pt}
\begin{tabularx}{\textwidth}{@{}>{\raggedright\arraybackslash}p{0.26\textwidth}
  >{\raggedright\arraybackslash}p{0.2\textwidth}
  >{\raggedright\arraybackslash}p{0.26\textwidth}
  >{\raggedright\arraybackslash}p{0.3\textwidth}@{}}
\toprule
Work & Range & Convexity Assumption& Quadratic Growth Assumption\\
\midrule
J\"orgens--Calabi--Pogorelov
\cite{Jorgens,Calabi,Pogorelov} & \(k=n\) & Convex & None \\
Bao--Chen--Guan--Ji \cite{BCGJ}
 & \(2\le k<n\) & \(D^2u>0\) & Lower \\
\addlinespace
Chang--Yuan \cite{CY}
 & \(k=2\) & Almost convex & None \\
\addlinespace
Li--Ren--Wang \cite{LRW}
 & \(2\le k<n\) & \((k+1)\)-convex & Lower \\
\addlinespace
Shankar--Yuan \cite{SY}
 & \(k=2\) & Semiconvex & None \\
\addlinespace
Wang--Bao \cite{WB}
 & \(1\le k\le n\) & \(k\)-admissible & Two-sided \\
\addlinespace
Chu--Dinew \cite{CD}
 & \(2\le k<n\) & \(k\)-admissible, \(\sigma_{k+1}\ge-A\) & Lower \\
\addlinespace
Tu \cite{Tu}
 & \(k=n-1\), \(n\ge3\) & \(k\)-admissible & Lower \\
\addlinespace
Zhang \cite{Zhang}
 & \(3\le k<n\) & Semiconvex & Lower \\
\addlinespace
B.~Wang \cite{BW}
 & \(2\le k<n\), \(2k>n\) & \(k\)-admissible & Lower \\
\midrule
\Cref{thm:main} & \(3\le k\le n-1\) & Convex & None \\
\bottomrule
\end{tabularx}
\end{table}

The distinction between convexity and admissibility is necessary.
Warren \cite{Warren} constructed nonpolynomial entire
\(k\)-admissible solutions when \(k\ge2\) and \(n\ge2k-1\).
Those solutions are nonconvex. \Cref{thm:main} shows that convex
solutions are quadratic for every \(n\) and \(k\).

\vspace{0.7\baselineskip}

We describe the proof. We map bounded sections
\(\{u<\ell+h\}\), with \(\ell\) affine and \(h>0\), between
fixed balls. The transformed equation is \(P(D^2w)=1\), where
\[
 P(H)=\sum_{|I|=k}\alpha_I\det H_I,
 \qquad \alpha_I\ge0,\quad\sum_I\alpha_I=1.
\]
Here \(H_I\) is a principal submatrix. Set
\(d_i=\sum_{I\ni i}\alpha_I\). We also include coefficient limits.
For smooth strictly convex solutions with zero boundary values,
we prove uniform interior Hessian bounds when \(d_i\ge\eta>0\).
We combine Legendre duality \cite{BCGJ}, a Jacobi inequality
for the inverse Hessian related to \cite{GS}, and integral
and contact set estimates \cite{Mooney25}.

Next we prove smoothness for convex viscosity solutions with
zero boundary values and \(d_i>0\). Trudinger's Dirichlet theory
\cite{Trudinger} supplies smooth convex approximations where every
supporting plane has only one contact point. Our estimates then
give regularity there.
Related regularity results appear in \cite{CJTZ,Hu}.
To exclude contact segments, we use sections as in \cite{CaffLoc}
and induct on dimension. In directions where \(d_i\to0\),
rescaling by \(\sqrt{d_i}\) makes the domains expand to whole lines.
The bounded convex limit is constant along these lines. Estimates of
Monge--Amp\`ere measures on slices preserve zero boundary values,
so the limit is smooth by induction. Savin's perturbation theorem
\cite{Savin} transfers regularity to nearby solutions.
The constant rank theorem \cite{BG} and zero boundary values
then exclude the contact segment.

Finally, remove the affine directions and choose large sections
with a fixed center. If some \(d_i\to0\), the limiting rescaled
Hessian vanishes in that direction at the center. The change of variables
instead gives a positive lower bound at the same point.
Thus all \(d_i\) stay uniformly positive, and the interior estimates
give a global Hessian bound.
By the Evans--Krylov estimate \cite{CS} on expanding balls,
\(D^2u\) is constant.

\Cref{sec:algebra} introduces the operators used in the proof.
We establish the classical estimates in \cref{sec:estimates}
and regularity for convex viscosity solutions in \cref{sec:weak}.
\Cref{sec:entire} completes the proof of \cref{thm:main}.

\section{Weighted operators and their elliptic branch}\label{sec:algebra}

We write \(H\ge0\) and \(H>0\) when the symmetric matrix \(H\)
is positive semidefinite and positive definite, respectively.
The identity matrix is \(\Id\), and \(\tr H\) is the sum of the
diagonal entries of \(H\). For symmetric matrices,
\(A:B=\sum_{i,j}A_{ij}B_{ij}=\tr(AB)\).
Thus \(DP(H):K\) is the derivative of \(P\) at \(H\) in direction
\(K\). We write \(\Sym_m\) for the real symmetric \(m\times m\)
matrices and \(B_r\) for the Euclidean ball of radius \(r\).

For an index set \(I\), the matrix \(H_I\) consists of the rows
and columns of \(H\) indexed by \(I\). The notation \(C\oplus0\)
means the block diagonal matrix with blocks \(C\) and zero.

\begin{definition}\label[definition]{def:class}
For \(1\le k\le m\), let \(\PP_{m,k}\) consist of the polynomials
\begin{equation}\label{eq:class}
 P_T(H)=\frac{\sigma_k(T^{1/2}HT^{1/2})}{\sigma_k(T)},
 \qquad T>0,
\end{equation}
and all their limits obtained by convergence of the coefficients.
In particular, \(P(\Id)=1\) for every \(P\in\PP_{m,k}\).
\end{definition}

After an orthogonal change of coordinates, every member has the form
\begin{equation}\label{eq:coeff}
 P(H)=\sum_{|I|=k}\alpha_I\det H_I,
\end{equation}
where \(\alpha_I\ge0\) and \(\sum_I\alpha_I=1\).
Indeed, diagonalize an approximating sequence \(T_j\) and pass to
a subsequence of its orthogonal matrices and of its normalized
minor coefficients. Set \(d_i=\sum_{I\ni i}\alpha_I\).
Then \(0\le d_i\le1\), \(\sum_i d_i=k\), and
\[
 DP(\Id)=\diag(d_1,\ldots,d_m).
\]

These polynomials are hyperbolic in direction \(\Id\): for every
symmetric \(H\), the polynomial \(t\mapsto P(H+t\Id)\) has only
real roots. Let \(\Gamma_P\) be the connected component of
\(\{P>0\}\) containing the positive definite matrices.
This is the hyperbolicity cone of \(P\). It is convex, its closure contains
the positive semidefinite cone, and \(P^{1/k}\) is concave on
\(\Gamma_P\). These properties hold for \cref{eq:class} and pass
to coefficient limits; see \cite{BGLS}. In particular,
\(H\ge0\) and \(P(H)>0\) imply \(H\in\Gamma_P\).
We will only use the class \(\PP_{m,k}\), or positive multiples
of its members, in statements involving these cone properties.

We will use two consequences of \cref{def:class}.
For \(P\in\PP_{m,k}\) and invertible \(L\),
\begin{equation}\label{eq:congruence}
 \frac{P(L^{-T}HL^{-1})}{P(L^{-T}L^{-1})}\in\PP_{m,k}.
\end{equation}
For an \(a\)-dimensional subspace \(V\), with \(a\ge k\) and
\(P(\Id_V\oplus0)>0\), we also have
\begin{equation}\label{eq:restriction}
 C\longmapsto
 \frac{P(C\oplus0_{V^\perp})}{P(\Id_V\oplus0_{V^\perp})}
 \ \in\PP_{a,k}.
\end{equation}
For \(P=P_T\), these operations replace \(T\) by
\(L^{-1}TL^{-T}\) and its compression to \(V\), respectively.
Both assertions pass to limits of coefficients, since their
denominators stay positive.

\begin{lemma}\label[lemma]{lem:logmono}
For \(P\in\PP_{m,k}\) and \(0<H_1\le H_2\),
\begin{equation}\label{eq:logmono}
 D\log P(H_1)\ge D\log P(H_2).
\end{equation}
Consequently, for \(H\ge0\) with \(P(H)>0\),
\begin{equation}\label{eq:ellipticity}
 DP(H)\ge\frac{P(H)}{\tr H}\,DP(\Id).
\end{equation}
For an unnormalized positive multiple the right side is
\(P(H)DP(\Id)/(P(\Id)\tr H)\).
\end{lemma}

\begin{proof}
A real polynomial is stable if it does not vanish when all
variables have positive imaginary parts \cite{BBL}.
First take \(P=P_T\). For \(H>0\) and vectors \(\xi,\zeta\), write
\[
 P(H+s\xi\xi^T+t\zeta\zeta^T)=a+bs+ct+dst.
\]
To check real stability, note that the determinant of
\[
 T^{1/2}(H+s\xi\xi^T+t\zeta\zeta^T)T^{1/2}+z\Id
\]
is nonzero when \(s,t,z\) have positive imaginary parts.
Differentiate \(m-k\) times in \(z\) and set \(z=0\).
By the Gauss--Lucas and Hurwitz theorems, the resulting polynomial
is stable; it is nonzero since \(a>0\).
Hence \(bc-ad\ge0\), since otherwise
\(s=-(a+ct)/(b+dt)\) lies in the upper half-plane when \(t\) does.
Thus
\[
 D^2\log P(H)[\xi\xi^T,\zeta\zeta^T]\le0.
\]
Decomposing a positive semidefinite increment into rank-one
matrices and integrating proves \cref{eq:logmono}.
Pass to limits of coefficients to obtain the general case.
By homogeneity and \(H\le(\tr H)\Id\), we obtain
\cref{eq:ellipticity} for \(H>0\), and then for \(H\ge0\)
by replacing \(H\) with \(H+\varepsilon\Id\).
\end{proof}

We need an operator defined also at indefinite test Hessians.
Let \(P\) be a positive multiple of a member of \(\PP_{m,k}\).
For \(M\in\Sym_m\), let \(t_P(M)\) be the unique real number such that
\[
 M+t_P(M)\Id\in\Gamma_P,\qquad
 P(M+t_P(M)\Id)=1,
\]
and set
\begin{equation}\label{eq:extension}
 E_P(M)=-t_P(M).
\end{equation}
This is well defined because \(P(M+t\Id)\) increases from zero
to infinity beyond its largest root. The function \(E_P\) is continuous,
nondecreasing under positive semidefinite increments, and
\(E_P(M+s\Id)=E_P(M)+s\).
In particular it is Lipschitz with respect to the operator norm.
If \(P_j\to P\) and \(P(\Id)>0\), then
\(E_{P_j}\to E_P\) locally uniformly.
This follows from continuity of polynomial roots.
By the implicit function theorem, at \(H\in\Gamma_P\) with \(P(H)=1\),
\begin{equation}\label{eq:DE}
 DE_P(H)=\frac{DP(H)}{\tr DP(H)}.
\end{equation}
In particular, if \(H\ge0\) and \(DP(\Id)>0\), this derivative
is positive definite by \cref{lem:logmono}.

\begin{definition}\label[definition]{def:viscosity}
A convex continuous function \(v\) is a viscosity solution of
\(P(D^2v)=1\) if \(E_P(D^2v)=0\) in the viscosity sense.
Thus a smooth function \(\phi\) touching \(v\) from above satisfies
\(E_P(D^2\phi)\ge0\) at the contact point.
A function touching from below satisfies
\(E_P(D^2\phi)\le0\).
The latter condition defines a viscosity supersolution.
Both conditions are imposed for all test Hessians, including
indefinite ones.
\end{definition}

At a positive semidefinite argument, the supersolution inequality
is equivalent to \(P(H)\le1\). After an invertible linear change
of variables, both viscosity inequalities hold for the transformed
operator. Their signs are
determined by the cone \(\Gamma_P\) and the surface \(P=1\),
although the value of \(E_P\) may change.
We use viscosity stability in this convention; see \cite{CIL}.
For comparison on a domain contained in \(B_R\), add
\(\delta(|x|^2-R^2)/2\) to the subsolution. The boundary ordering
is preserved and the subsolution inequality becomes strict.
A positive interior maximum would give viscosity jets \(X\le Y\)
with \(\delta\le E_P(X)\le E_P(Y)\le0\), a contradiction.
Let \(\delta\downarrow0\).

\section{Classical interior estimates}\label{sec:estimates}

We first estimate smooth solutions on a normalized convex domain.
The constants must remain uniform as the coefficients of \(P\)
vary, provided \(DP(\Id)\) stays positive definite.

\begin{theorem}\label{thm:classical}
Let \(2\le k<m\), \(P\in\PP_{m,k}\), and
\(DP(\Id)\ge\eta\Id\) for some \(\eta>0\).
Let \(\Omega\) be an open convex set with
\(B_1\subset\Omega\subset B_R\).
Suppose \(w\in C^4(\Omega)\cap C(\overline\Omega)\) satisfies
\(D^2w>0\) and
\[
 \begin{cases}
 P(D^2w)=1 & \text{in }\Omega,\\
 w=0 & \text{on }\partial\Omega.
 \end{cases}
\]
There exists \(C=C(m,k,R,\eta)\) such that
\begin{equation}\label{eq:classicalC2}
 C^{-1}\Id\le D^2w\le C\Id\qquad\text{on }B_{1/2}.
\end{equation}
\end{theorem}

We first bound \(\tr(D^2w)^{-1}=\sum_i1/\lambda_i\), where
\(\lambda_i>0\) are the Hessian eigenvalues. This bounds them from
below; the equation then bounds them from above.
We use a Jacobi inequality, integral bounds, and a contact set argument.

\subsection{A Jacobi inequality for the inverse Hessian}

We use Legendre duality as in \cite{BCGJ} and a calculation
related to the concavity property in \cite{GS}. The calculation
allows zero coefficients.

\begin{lemma}\label[lemma]{lem:jacobi}
Let \(P(H)=\sum_{|I|=k}\omega_I\det H_I\), where
\(\omega_I\ge0\) and at least one coefficient is positive.
Suppose \(w\in C^4\), \(H=D^2w>0\), and \(P(H)=1\).
Set \(F=DP(H)\) and \(U=\tr H^{-1}\). For a smooth function \(f\),
write \(Lf=\sum_{i,j}F^{ij}f_{ij}\) and
\(\lvert Df\rvert_F^2=\sum_{i,j}F^{ij}f_i f_j\). Then
\begin{equation}\label{eq:jacobi}
 L\log U\ge\frac1k\abs{D\log U}_F^2.
\end{equation}
Moreover \(M\mapsto\log P(M^{-1})\) is convex on \(M>0\).
\end{lemma}

\begin{proof}
Take the local Legendre transform \(W=w^*\), and put
\(M=D^2W=H^{-1}\) and \(r=m-k\).
Set \(c_J=\omega_{J^c}\) and
\[
 Q(M)=\sum_{|J|=r}c_J\det M_J,
\]
with the empty determinant equal to one.
By the identity for complementary principal minors,
\(P(M^{-1})=Q(M)/\det M\). Set \(G=\log\det M-\log Q(M)\).
For the terms with \(c_J>0\), define
\[
 \pi_J=\frac{c_J\det M_J}{Q(M)},\qquad
 \Pi_J=M^{1/2}E_JM_J^{-1}E_J^TM^{1/2},
\]
where the columns of \(E_J\) are the coordinate vectors indexed
by \(J\). The weights satisfy \(\pi_J\ge0\) and
\(\sum_J\pi_J=1\). Each \(\Pi_J\) is an orthogonal projection
of rank \(r\). Write
\(B=\sum_J\pi_J\Pi_J\), \(R_J=\Id-\Pi_J\).
For \(Z=M^{-1/2}XM^{-1/2}\), we obtain
\[
 A:=DG(M)=M^{-1/2}(\Id-B)M^{-1/2}.
\]
The second derivative is
\begin{equation}\label{eq:variance}
 \begin{aligned}
 -D^2G(M)[X,X]
 &=\tr Z^2-\sum_J\pi_J\tr(\Pi_JZ\Pi_JZ)\\
 &\quad+\sum_J\pi_J\bigl(\tr(\Pi_JZ)-\tr(BZ)\bigr)^2.
 \end{aligned}
\end{equation}
Since \(\Pi_J\) is an orthogonal projection,
\(\tr(\Pi_JZ\Pi_JZ)\le\tr Z^2\). Thus the right side of
\cref{eq:variance} is nonnegative, and \(-G\) is convex.

Suppose now \(DG(M)[X]=0\). Then \(\tr(BZ)=\tr Z\), so the
last sum in \cref{eq:variance} is
\(\sum_J\pi_J(\tr(R_JZ))^2\).
In a suitable orthonormal basis,
\[
 Z=\begin{pmatrix}C&D\\D^T&E\end{pmatrix},
 \qquad R_J=\begin{pmatrix}\Id_k&0\\0&0\end{pmatrix}.
\]
The corresponding contribution to \cref{eq:variance} is
\[
 a_J=\tr C^2+2\tr(DD^T)+(\tr C)^2.
\]
For a symmetric \(k\times k\) matrix \(C\), let \(\lambda\) be
an eigenvalue of largest absolute value. Then
\[
 \tr C^2+(\tr C)^2\ge\frac{k+1}{k}\lambda^2.
\]
For \(k>1\), minimize
\(\lambda^2+s^2/(k-1)+(\lambda+s)^2\), where \(s\) is the sum
of the other eigenvalues. The case \(k=1\) is immediate.
The nonzero eigenvalues of \(ZR_JZ\) are those of \(C^2+DD^T\).
As \((k+1)/k\le2\), this proves
\[
 a_J\Id\ge(1+1/k)ZR_JZ.
\]
Multiply by \(\pi_J\) and sum over \(J\). We obtain
\begin{equation}\label{eq:strongmatrix}
 [-D^2G(M)[X,X]]\Id\ge(1+1/k)Z(\Id-B)Z
 \quad\text{if }DG(M)[X]=0.
\end{equation}

For a fixed unit vector \(e\), let \(h=W_{ee}>0\) and
\(X=D^2W_e\). By differentiating \(G(D^2W)=0\) twice, we obtain
\[
 DG(M)[X]=0,\qquad
 A:D^2h=-D^2G(M)[X,X].
\]
By symmetry of the third derivatives, \(Dh=Xe\).
Test \cref{eq:strongmatrix} on the unit vector \(M^{1/2}e/\sqrt h\).
It follows that
\[
 A:D^2h\ge(1+1/k)\frac{\abs{Dh}_{A}^2}{h}.
\]
Sum over an orthonormal basis and use Cauchy--Schwarz to obtain
the same inequality for \(\tr M\). Under \(y=Dw(x)\),
\(A=HFH\) and \(Lw_i=0\), so there is no drift term.
Taking the logarithm and changing back to \(x\) proves
\cref{eq:jacobi}.
\end{proof}

\subsection{Integral bounds}

We next bound the Hessian integrals and
\(\log\tr(D^2w)^{-1}\), uniformly when \(DP(\Id)\ge\eta\Id\).

\begin{lemma}\label[lemma]{lem:integrals}
Let \(2\le k<m\), \(P\in\PP_{m,k}\), \(DP(\Id)\ge\eta\Id\),
and \(B_1\subset\Omega\subset B_R\), with \(\Omega\) convex.
Suppose \(w\in C^4(\Omega)\cap C(\overline\Omega)\),
\(D^2w>0\), \(P(D^2w)=1\), and \(w=0\) on \(\partial\Omega\).
For each fixed \(0<s<1/2\), set \(D_s=\{w<-s\}\).
There are constants depending only on \(m,k,R,\eta,s\) such that
\begin{gather}
 -R^2/2\le w\le0,\qquad w(0)\le-1/2,\label{eq:C0}\\
 \dist(D_s,\partial\Omega)\ge c,\qquad
 \sup_{D_s}|Dw|\le C,\label{eq:depth}\\
 \int_{D_s}\tr D^2w\dd\le C,\qquad
 \int_{D_s}\det D^2w\dd\ge c,\label{eq:mass}\\
 \int_{D_s}(\det D^2w)^{1+1/(m-k)}\dd\le C,\label{eq:massLp}\\
 \int_{D_s}\log\tr(D^2w)^{-1}\dd\le C.\label{eq:logL1}
\end{gather}
\end{lemma}

\begin{proof}
Quadratic comparison proves \cref{eq:C0}.
We first choose coordinates in which
\begin{equation}\label{eq:star}
 P(H)\ge c_*\sum_{j\in E}\det H_{A\cup\{j\}},
 \qquad H\ge0,
\end{equation}
where \(A=\{1,\ldots,k-1\}\), \(E=A^c\), and
\(c_*=c_*(m,k,\eta)>0\).
For \(P=P_T\), order the eigenvalues \(t_i\) of \(T\) decreasingly
and put \(N=\binom mk\). Then
\[
 \eta\le d_m
 =\frac{t_m\sigma_{k-1}(t\mid m)}{\sigma_k(t)}
 \le N\frac{t_m}{t_k}.
\]
Since \(\sigma_k(t)\le Nt_1\cdots t_k\), each coefficient in
\cref{eq:star} is at least \(\eta/N^2\).
For a limit polynomial, take approximants with
\(DP_T(\Id)\ge(\eta/2)\Id\) and pass to a subsequence of their
orthogonal coordinates. Thus \(c_*=\eta/(2N^2)\) suffices.

Put \(r=m-k\). We will also use
\begin{equation}\label{eq:detLp}
 (\det H)^{1+1/r}\le C\det H_E,
 \qquad \text{for }H>0,\ P(H)=1.
\end{equation}
Indeed, for \(d=\det H\) and \(M=H^{-1}\), complementary
minors and Maclaurin's inequality give
\[
 1\ge c_*d\,\sigma_r(M_E)
 \ge c_*d\,\sigma_r(H_E^{-1})
 \ge c_*(r+1)d\,(\det H_E)^{-r/(r+1)},
\]
where \(M_E\ge H_E^{-1}\).

To prove \cref{eq:depth}, note that the polynomial
\[
 P_*(H)=\sum_{j\in E}\det H_{A\cup\{j\}}
       =\det H_A\,\tr(H/H_A)
\]
is invariant under separate orthogonal transformations of \(A,E\);
\(H/H_A\) denotes the Schur complement.
For every unit vector \(\nu\), choose a \(k\)-plane \(V\) that
contains \(A\) and \(\nu\). By \cref{eq:star},
\(P(H)\ge c_*\det(H|_V)\). Choose \(\nu\) normal to a nearest supporting plane.
The restriction to \(\Omega\cap(x+V)\) has zero boundary values,
diameter at most \(2R\), and Monge--Amp\`ere density at most
\(c_*^{-1}\). By the Alexandrov maximum principle,
\[
 |w(x)|^k\le C\dist(x,\partial\Omega).
\]
This proves \cref{eq:depth}, including the gradient bound on a
fixed neighborhood of \(D_s\) by convexity.

For the Hessian integrals, we use Newton tensors as in \cite{CW}.
The Newton tensor \(T_{j-1}=D\sigma_j(D^2w)\) is positive
semidefinite and divergence free, and
\[
 T_{j-1}:D^2w=j\sigma_j(D^2w),\qquad
 \tr T_{j-1}=(m-j+1)\sigma_{j-1}(D^2w).
\]
For a nonnegative compactly supported cutoff \(\chi\),
\[
 \int\chi\sigma_j(D^2w)
 =-\frac1j\int D\chi\cdot T_{j-1}Dw.
\]
Using nested cutoffs, we obtain
\(\int_{D_s}\sigma_j(D^2w)\le C\) for every \(j\).
For sufficiently small \(|p|\), the minimum of \(w-p\cdot x\)
on \(\overline D_s\) is interior: its value at zero is at most
\(-1/2\), whereas its boundary value is at least \(-s-R|p|\).
Thus \(Dw(D_s)\) contains a fixed ball. Strict convexity and the
area formula prove the second assertion in \cref{eq:mass}.

The partial gradient on each \(E\)-parallel slice of \(D_s\)
is injective and bounded by \cref{eq:depth}. Integrate
\cref{eq:detLp} and apply the area formula on these slices.
Since the projection onto \(A\) is bounded, we obtain \cref{eq:massLp}.

Put \(H=D^2w\), \(F=DP(H)\), \(S=\tr H\), and
\(g=\log\tr H^{-1}\). We have \(g\ge0\), since \(H>\Id\) would
contradict \(P(H)=P(\Id)=1\).
The matrix \(F\) is divergence free, that is,
\(\sum_i\partial_i F^{ij}=0\) for every \(j\).
This follows from cancellation of mixed derivatives in the cofactors.
Multiply \cref{eq:jacobi} by \(\chi^2\) and integrate by parts.
We obtain
\[
 \frac1k\int\chi^2|Dg|_F^2
 \le-2\int\chi\, D\chi\cdot F Dg
 \le\frac1{2k}\int\chi^2|Dg|_F^2
       +2k\int|D\chi|_F^2.
\]
In the coordinates \cref{eq:coeff},
\(\tr F\le\sigma_{k-1}(H)\), since each coefficient of a
\((k-1)\)-minor in \(\tr F\) is a subsum of the \(\alpha_I\).
Thus \(\int_{D_s}|Dg|_F^2\le C\).
By \cref{eq:ellipticity},
\[
 \int_{D_s}|Dg|
 \le \eta^{-1/2}
       \left(\int_{D_s}|Dg|_F^2\right)^{1/2}
       \left(\int_{D_s}S\right)^{1/2}\le C.
\]

By \cref{eq:mass,eq:massLp},
\(\det H\ge c_1>0\) on a set of volume at least \(c_2>0\).
Choose \(K\) using the \(L^1\) bound for \(S\), so that removing
\(\{S>K\}\) leaves a set \(E_0\) of volume at least \(c_2/2\).
On \(E_0\),
\(\tr H^{-1}=\sigma_{m-1}(H)/\det H\le C\), so \(g\le C\).
Let \(\bar g\) be the average of \(g\)
over \(D_s\). By the \(L^1\) Poincar\'e inequality on a convex
domain \cite{AD} and the bound for \(\int_{D_s}|Dg|\),
\[
 \bar g|E_0|\le\int_{E_0}g+
            \int_{D_s}|g-\bar g|\le C|E_0|+C.
\]
Together with \(|E_0|\ge c_2/2\) and \(|D_s|\le|B_R|\),
this proves \cref{eq:logL1}.
\end{proof}

\subsection{The pointwise Hessian bound}

We now use a double logarithm and the Alexandrov--Bakelman--Pucci
method to pass from integral to pointwise bounds; see \cite{Mooney25}
for related arguments for \(\sigma_2\).

\begin{proof}[Proof of \cref{thm:classical}]
Put \(H=D^2w\), \(F=DP(H)\), \(N=\tr F\), and \(S=\tr H\).
Hyperbolicity and normalization give
\[
 P(H+t\Id)=\prod_{\ell=1}^k(t+\mu_\ell),\quad
 \mu_\ell>0,\quad \prod_\ell\mu_\ell=1,\quad
 N=\sum_\ell\mu_\ell^{-1}.
\]
The largest hyperbolic eigenvalue is nondecreasing in the
semidefinite order. Indeed,
\[
 \mu_{\max}(H)=\inf\{c:c\Id-H\in\overline{\Gamma_P}\}.
\]
If \(\lambda=\lambda_{\max}(H)\) has unit eigenvector \(e\),
then \(H\ge\lambda ee^T\). Expanding along the rank-one matrix
\(ee^T\), we obtain
\[
 P(t\Id+\lambda ee^T)
 =t^{k-1}\big(t+\lambda DP(\Id)[ee^T]\big).
\]
Consequently \(\mu_{\max}(H)\ge\eta\lambda\).
By the arithmetic--geometric mean inequality for the other
\(k-1\) inverse roots,
\begin{equation}\label{eq:tracefromN}
 S\le\frac m\eta\left(\frac N{k-1}\right)^{k-1}.
\end{equation}

On \(D=\{w<-1/4\}\), put
\[
 \rho=-w-1/4,\qquad g=\log\tr H^{-1},\qquad q=1+g.
\]
Then \(q\ge1\), \(\rho(0)\ge1/4\), \(0<\rho\le R^2/2\),
and \(\int_Dq\le C\) by \cref{lem:integrals}.
By the chain rule and \cref{eq:jacobi},
\begin{equation}\label{eq:loglog}
 L\log q\ge(q/k-1)|D\log q|_F^2.
\end{equation}
Set \(\beta=mk\), \(a=1/[4(1+2R)^2]\), and \(\delta=a\).
Consider the auxiliary function
\[
 \mathcal B=\log q+\beta\log\rho+a|x|^2.
\]
For each fixed solution, \(q\) is bounded on \(\overline D\);
thus \(\mathcal B\to-\infty\) at its boundary.
For every \(|p|\le\delta\), take all maxima of
\(\mathcal B-p\cdot x\). Their union \(\mathcal C\) is compact
in \(D\). At every such maximum, \(D\mathcal B=p\) and
\(D^2\mathcal B\le0\). By comparison with the value at zero, we obtain
\[
 q(x)\ge c_0q(0)\rho(x)^{-\beta},
\]
where \(c_0=4^{-\beta}e^{-aR^2-\delta R}\).
We may assume \(q(0)\) is large enough that \(q\ge2k\) on
\(\mathcal C\), since otherwise the desired bound already holds.

At a contact point,
\(D\log q=p-2ax-\beta D\rho/\rho\).
Using \(L\rho=-k\) and \cref{eq:loglog}, we obtain
\begin{align*}
 L\mathcal B
 &\ge(\beta^2/2-\beta)\frac{|D\rho|_F^2}{\rho^2}
       +\big(2a-(\delta+2aR)^2\big)N-\frac{\beta k}{\rho}\\
 &\ge aN-\frac{\beta k}{\rho}.
\end{align*}
Since \(L\mathcal B\le0\), \cref{eq:tracefromN} implies
\[
 N\le C/\rho,\qquad S\le C\rho^{-(k-1)}
 \quad\text{on }\mathcal C.
\]
By \cref{eq:ellipticity},
\(\det F\ge c\rho^{m(k-1)}\).
Set \(K=-D^2\mathcal B\ge0\) at these points. Then
\[
 \det K\le\frac{(\tr(FK)/m)^m}{\det F}
       \le C\rho^{-mk}=C\rho^{-\beta}.
\]
The gradient \(D\mathcal B(\mathcal C)\) contains \(B_\delta\).
By the area formula,
\[
 |B_\delta|\le\int_{\mathcal C}|\det D^2\mathcal B|
       \le C\int_{\mathcal C}\rho^{-\beta}
       \le\frac C{q(0)}\int_Dq\le\frac C{q(0)}.
\]
It follows that \(q(0)\le C\), hence \(H(0)\ge b\Id\) for
a uniform \(b>0\). For a largest eigenvector \(e\),
\[
 1=P(H)\ge
 P\big(b\Id+(\lambda_{\max}(H)-b)ee^T\big)
 \ge b^{k-1}(\lambda_{\max}(H)-b)\eta.
\]
This also bounds \(H(0)\) above.
For \(x_0\in B_{1/2}\), apply the estimate at zero to
\(\widetilde w(y)=4w(x_0+y/2)\). Its Hessian is unchanged and
its domain lies between \(B_1\) and \(B_{2R+1}\).
This proves \cref{eq:classicalC2}.
\end{proof}

\section{Regularity of convex viscosity solutions}\label{sec:weak}

\begin{theorem}\label{thm:weak}
Let \(1\le k\le m\), \(P\in\PP_{m,k}\), and \(DP(\Id)>0\).
Let \(\Omega\) be a bounded open convex set.
Every convex \(v\in C(\overline\Omega)\) satisfying
\[
 \begin{cases}
 P(D^2v)=1 & \text{in }\Omega,\\
 v=0 & \text{on }\partial\Omega
 \end{cases}
\]
in the viscosity sense of \cref{def:viscosity} is smooth in
\(\Omega\) and satisfies \(D^2v>0\).
\end{theorem}

We first prove regularity where supporting planes have only one
contact point. We then exclude nontrivial contact sets by induction
on dimension.

\subsection{Regularity at points of strict contact}\label{sec:strictcontact}

We construct convex approximations using the Legendre transform
\cite{BCGJ} and the Dirichlet theory for Hessian quotients
\cite{Trudinger}, then apply \cref{thm:classical}.

\begin{proposition}\label[proposition]{prop:strictcontact}
Let \(2\le k<m\), \(P\in\PP_{m,k}\), and \(DP(\Id)>0\).
Let \(v\in C(\overline\Omega)\) be a convex viscosity
solution of \(P(D^2v)=1\) on a bounded convex domain.
If every supporting affine function at \(x_0\in\Omega\)
has contact set \(\{x_0\}\) in \(\overline\Omega\), then \(v\)
is smooth and has positive definite Hessian near \(x_0\).
\end{proposition}

\begin{proof}
Choose \(P_j=P_{T_j}\) with \(T_j>0\) so that their coefficients
converge to those of \(P\). We may assume
\(DP_j(\Id)\ge\eta\Id/2\), where \(DP(\Id)\ge\eta\Id\).
By definition, \(P_j(\Id)=1\). Define, for \(M>0\),
\[
 Q_j(M)=P_j(M^{-1})^{-1/k},\qquad
 Q(M)=P(M^{-1})^{-1/k}.
\]
Extend both functions by zero outside the positive definite cone.
If \(0<M\le L\Id\) and
\(\lambda=\lambda_{\min}(M)\le L/2\), we claim that
\begin{equation}\label{eq:dualboundary}
 P_j(M^{-1})\ge\frac{\eta}{4}L^{-(k-1)}\lambda^{-1}.
\end{equation}
Indeed \(M^{-1}\ge L^{-1}\Id+
(\lambda^{-1}-L^{-1})ee^T\) for a unit minimum eigenvector \(e\).
Thus the zero extensions are continuous, and
\(Q_j\to Q\) locally uniformly on all symmetric
matrices.

For \(P_j=P_{T_j}\), put \(Z_j=\sigma_k(T_j)\) and
\(N=T_j^{-1/2}MT_j^{-1/2}\). Then
\begin{equation}\label{eq:quotient}
 Q_j(M)^k
 =Z_j\frac{\det N}{\sigma_{m-k}(N)}.
\end{equation}
The Hessian quotient root is concave on the positive definite
cone \cite{Trudinger}. Hence the \(Q_j\), and their
limit, are concave and homogeneous of degree one there, with
\(Q_j(\Id)=Q(\Id)=1\). In particular
\begin{equation}\label{eq:dualstrict}
 Q(M+t\Id)\ge Q(M)+t,
 \quad M\ge0,\ t\ge0.
\end{equation}
For comparison, add \(\delta(|p|^2-R^2)/2\) to a convex
subsolution on a domain contained in \(B_R\).
Its upper viscosity jets satisfy \(X\ge\delta\Id\), and by
\cref{eq:dualstrict}, \(Q(X)\ge1+\delta\).
At a positive interior maximum, \(X\le Y\) and \(Q(Y)\le1\),
a contradiction. Continuity covers closed jets.
Letting \(\delta\downarrow0\) proves comparison for \(Q\), and
likewise for each \(Q_j\).

Extend \(v\) by \(+\infty\) outside \(\overline\Omega\).
Its Legendre transform is
\[
 U(p)=\max_{x\in\overline\Omega}\{p\cdot x-v(x)\}.
\]
To localize the dual equation, set
\(b(p)=\max_{x\in\partial\Omega}\{p\cdot x-v(x)\}\)
and \(\mathcal O=\{U>b\}\).
The function \(U\) is finite, convex, and globally Lipschitz.
All maximizers for \(p\in\mathcal O\) lie in \(\Omega\).
We claim that \(Q(D^2U)=1\) on the open set \(\mathcal O\).
If \(\psi\) touches \(U\) from above at \(p_0\), convexity makes
\(U\) differentiable there, with \(x_0=D\psi(p_0)\in\Omega\),
and \(M=D^2\psi(p_0)\ge0\). Add
\(\delta|p-p_0|^2/2\) and take its local strictly convex conjugate.
That conjugate touches \(v\) from below at \(x_0\) and has
Hessian \((M+\delta\Id)^{-1}\). The supersolution condition gives
\(Q(M+\delta\Id)\ge1\); let \(\delta\downarrow0\).
In particular \cref{eq:dualboundary} excludes a singular upper test.

If \(\psi\) touches \(U\) from below and \(D^2\psi(p_0)>0\),
then \(x_0=D\psi(p_0)\in\partial U(p_0)\subset\Omega\).
The strict quadratic lower test and convexity make \(p_0\) the
unique minimizer of \(U(p)-x_0\cdot p\).
This function grows at least linearly at infinity because
\(x_0\in\Omega\). Hence conjugate maximizers for \(x\) near
\(x_0\) stay near \(p_0\). The local conjugate of \(\psi\)
touches \(v\) from above, so \(Q(D^2\psi(p_0))\le1\).
For lower tests whose Hessians are not positive definite,
the inequality follows from the zero extension. The claim is proved.

The contact hypothesis implies that the compact convex set
\(K=\partial v(x_0)\) lies in \(\mathcal O\).
Choose a smooth uniformly convex domain
\(K\Subset D\Subset\mathcal O\) and smooth boundary data
\(\psi_a\to U|_{\partial D}\) uniformly.
For every fixed \(j,a\), \cref{eq:quotient} and the Dirichlet
theorem \cite{Trudinger} give a smooth strictly
convex solution
\[
 \begin{cases}
 Q_j(D^2U_{j,a})=1 & \text{in }D,\\
 U_{j,a}=\psi_a & \text{on }\partial D.
 \end{cases}
\]
The admissible cone for this quotient equation is \(\Gamma_m\),
so \(D^2U_{j,a}>0\).
The fixed affine change associated with \(T_j\) preserves
uniform convexity of \(D\).

For fixed \(a\), compactness up to \(\partial D\) follows from
barriers independent of \(j\). With \(\nu_y\) the outer normal,
uniform convexity gives
\(\nu_y\cdot(p-y)\le-\kappa|p-y|^2\) on \(\partial D\).
For \(B\ge1\) fixed and \(A\) large, the functions
\begin{align*}
 q_y(p)&=\psi_a(y)+D\psi_a(y)\cdot(p-y)
                 +A\nu_y\cdot(p-y)+B|p-y|^2,\\
 \ell_y(p)&=\psi_a(y)+D\psi_a(y)\cdot(p-y)
                 -A\nu_y\cdot(p-y)
\end{align*}
bracket \(\psi_a\) on \(\partial D\) and have operator values
\(2B\) and zero. They give uniform bounds and boundary continuity;
convexity gives interior equicontinuity.
By stability and comparison, \(U_{j,a}\to U_a\) uniformly and
\[
 \norm{U_a-U}_{C^0(\overline D)}
 \le\norm{\psi_a-U}_{C^0(\partial D)}.
\]
Choose a diagonal sequence of smooth strictly convex solutions
\(U_j\) of \(Q_j(D^2U_j)=1\) such that
\(U_j\to U\) uniformly on \(\overline D\).

Since the minimum set of \(U(p)-x_0\cdot p\) is \(K\Subset D\),
its boundary minimum is strictly larger. This gap persists for
\(x\) near \(x_0\) and large \(j\). Thus
\[
 v_j(x)=\max_{p\in\overline D}\{x\cdot p-U_j(p)\}
\]
is smooth and strictly convex on a fixed ball \(B_r(x_0)\),
satisfies \(P_j(D^2v_j)=1\), and converges uniformly there to \(v\).
Indeed, the same gap for \(U\) places its minimizers in \(D\).
By convexity they are global, so the restricted Legendre transform
recovers \(v\).

Choose any supporting affine function \(\ell\) at \(x_0\).
Its unique contact gives
\(\min_{\partial B_r(x_0)}(v-\ell)>0\).
For a fixed small \(h>0\), the sections
\[
 B_\rho(x_0)\subset\Omega_j=\{v_j-\ell<h\}
                    \Subset B_r(x_0)
\]
therefore have a fixed \(\rho>0\) for all large \(j\).
The functions \(v_j-\ell-h\) solve \(P_j(D^2w)=1\) classically
on \(\Omega_j\), with positive definite Hessian and zero boundary
values. The rescaling
\(\rho^{-2}(v_j(x_0+\rho y)-\ell(x_0+\rho y)-h)\)
preserves the Hessian and the equation.
\Cref{thm:classical} gives uniform two-sided Hessian bounds near
\(x_0\). Concavity, uniform ellipticity, the Evans--Krylov estimate
\cite{CS}, and successive Schauder estimates imply smooth convergence
on a smaller ball. This proves the proposition.
\end{proof}

By \cref{prop:strictcontact}, it remains to exclude a segment of
contact points. We normalize sections near an exposed point of
the contact set. If the coefficients degenerate, a further rescaling
produces a solution in fewer variables. Its zero boundary values
will allow induction on dimension.

\subsection{Boundary values under degeneration of the coefficients}

The Monge--Amp\`ere measure of a convex function \(f\) is
\(\MA[f](E)=|\partial f(E)|\) for Borel sets \(E\).
We use its weak continuity under local uniform convergence
and the Alexandrov maximum principle; see
\cite{SavinObstacle}.

\begin{lemma}\label[lemma]{lem:slices}
Let \(P(H)=\sum_{|I|=k}\alpha_I\det H_I\) be a positive
multiple of a member of \(\PP_{m,k}\).
Let \(v\) be finite and convex on an open convex set \(\Omega\),
and a supersolution in the sense of \cref{def:viscosity}.
For every \(\alpha_I>0\) and every affine plane parallel to the
coordinate \(I\)-plane,
\begin{equation}\label{eq:slicemeasure}
 \MA[v|_{a+\R^I}]\le\alpha_I^{-1}\,dx_I
 \quad\text{in }\Omega\cap(a+\R^I).
\end{equation}

If in addition \(\Omega\) has diameter at most \(2R\),
\(v\in C(\overline\Omega)\), and \(v=0\) on \(\partial\Omega\),
the following boundary bounds hold.
Let \(e\) be a unit normal to a supporting hyperplane, with
\(e^TDP(\Id)e>0\). If \(t\) is the distance of \(x\in\Omega\)
to that hyperplane, then
\[
 |v(x)|^k\le C(m,k,R)\frac{t}{e^TDP(\Id)e}.
\]
The same argument applies to a bounded coordinate slice with
index set \(A\), using only coefficients with \(I\subset A\).
In particular, if
\(\sum_{I\subset A,\ i\in I}\alpha_I\ge\eta>0\) for every \(i\in A\),
then on such a slice
\begin{equation}\label{eq:activeboundary}
 |v(x)|^k\le C(m,k,R,\eta)\,
              \dist(x,\partial(\Omega\cap(a+\R^A))).
\end{equation}
Here only the diameter of the \(A\)-slice needs to be bounded.
\end{lemma}

\begin{proof}
Locally, extend \(v\) from a larger bounded convex subdomain
by \(+\infty\) and take its Moreau infimal convolution
\(v_\varepsilon\) \cite{Moreau}. It is convex and \(C^{1,1}\),
converges locally uniformly to \(v\), and has minimizing points
inside \(\Omega\) for small \(\varepsilon\).
If a smooth function \(\phi\) touches \(v_\varepsilon\) from below
at \(x\), and \(y\) is the minimizing point, then
\[
 z\longmapsto\phi(z+x-y)-\frac{|x-y|^2}{2\varepsilon}
\]
touches \(v\) from below at \(y\). Thus \(v_\varepsilon\) retains
the supersolution condition on the smaller set.
At a point of second differentiability, lower quadratic tests
with Hessian \(D^2v_\varepsilon-\delta\Id\), followed by
\(\delta\downarrow0\), give
\(P(D^2v_\varepsilon)\le1\).
This includes singular positive semidefinite Hessians.
Nonnegativity of all principal minors implies
\[
 \alpha_I\det D_I^2v_\varepsilon\le1
 \quad a.e..
\]

By Fubini's theorem and the area formula for the Lipschitz gradient,
\cref{eq:slicemeasure} holds for \(v_\varepsilon\) on almost every
parallel slice. Approximate an exceptional slice by such slices.
The restrictions converge uniformly on compact subsets, so weak
continuity of Monge--Amp\`ere measures extends the bound to every
slice. Let \(\varepsilon\downarrow0\) and exhaust the slice to
obtain \cref{eq:slicemeasure} for \(v\).

For the boundary estimate, let \(\pi_I\) denote orthogonal
projection onto the \(I\)-plane. Since
\[
 e^TDP(\Id)e=\sum_{|I|=k}\alpha_I|\pi_Ie|^2,
\]
some \(I\) satisfies
\(\alpha_I|\pi_Ie|^2\ge e^TDP(\Id)e/\binom mk\).
The restriction of \(v\) to \(\Omega\cap(x+\R^I)\) has zero
boundary values: every point on its relative boundary belongs
to \(\partial\Omega\).
Moving from \(x\) in the projected normal direction reaches the
supporting hyperplane after distance \(t/|\pi_Ie|\), and exits
the slice no later. By \cref{eq:slicemeasure} and the
Alexandrov maximum principle in dimension \(k\),
\[
 |v(x)|^k
 \le C R^{k-1}\frac{t}{|\pi_Ie|}
          \MA[v|_{x+\R^I}](\Omega\cap(x+\R^I))
 \le C\frac{t}{\alpha_I|\pi_Ie|}
 \le C\frac{t}{e^TDP(\Id)e}.
\]
For an \(A\)-slice, use only \(I\subset A\).
The denominator is at least \(\eta\); a nearest supporting
hyperplane then proves \cref{eq:activeboundary}.
\end{proof}

Rescaling by the square roots of the diagonal entries of
\(DP_j(\Id)\) stretches the directions whose entries tend to zero.
The limit is constant in those directions, while \cref{lem:slices}
preserves its boundary values in the remaining variables.
We will also use this limit for entire solutions.

\begin{lemma}\label[lemma]{lem:newton}
Let \(P_j\in\PP_{m,k}\) with \(DP_j(\Id)>0\).
Let \(\Omega_j\) be open convex sets satisfying
\(B_1\subset\Omega_j\subset B_R\), and let
\(w_j\in C(\overline\Omega_j)\) be convex viscosity solutions of
\[
 \begin{cases}
 P_j(D^2w_j)=1 & \text{in }\Omega_j,\\
 w_j=0 & \text{on }\partial\Omega_j.
 \end{cases}
\]
In orthogonal coordinates for \cref{eq:coeff}, write
\[
 P_j(H)=\sum_{|I|=k}\alpha_{I,j}\det H_I,\quad
 d_{i,j}=\sum_{I\ni i}\alpha_{I,j},\quad
 D_j=\diag(d_{i,j}).
\]
Suppose \(\min_i d_{i,j}\to0\). Pass to a subsequence such that
\(d_{i,j}\to d_i^*\), and set
\(A=\{i:d_i^*>0\}\), \(B=A^c\), and \(a=|A|\).
Then \(k\le a<m\). The functions
\[
 W_j(z)=w_j(D_j^{1/2}z),\qquad
 \mathcal D_j=D_j^{-1/2}\Omega_j
\]
converge locally uniformly, on a limiting cylinder
\(\Omega_*\times\R^B\), to \(\varphi(z_A)\).
Here \(\Omega_*\subset\R^A\) is bounded and convex,
\(B_1^a\subset\Omega_*\), and
\(\varphi\in C(\overline\Omega_*)\) is convex, zero on its boundary,
and \(\varphi(0)\le-1/2\).
Define the limiting coefficients by
\[
 \beta_{I,j}=\frac{\alpha_{I,j}}{\prod_{i\in I}d_{i,j}}
       \longrightarrow\beta_I^*\in[0,1].
\]
Set
\[
 \widehat Q(C)=\sum_{\substack{I\subset A\\|I|=k}}\beta_I^*\det C_I.
\]
Then \(\varphi\) solves
\begin{equation}\label{eq:activePDE}
 \widehat Q(D^2\varphi)=1\qquad\text{in }\Omega_*.
\end{equation}
Moreover \(\widehat Q(\Id)>0\),
\(D\widehat Q(\Id)\ge\Id\), and
\(\widehat Q/\widehat Q(\Id)\in\PP_{a,k}\).
The full limiting polynomial
\(Q_\infty(H)=\sum_I\beta_I^*\det H_I\) satisfies
\begin{equation}\label{eq:fullcoverage}
 1\le Q_\infty(\Id)\le\binom mk,\qquad DQ_\infty(\Id)\ge\Id.
\end{equation}
\end{lemma}

\begin{proof}
The transformed equation is \(Q_j(D^2W_j)=1\), where
\[
 Q_j(H)=P_j(D_j^{-1/2}HD_j^{-1/2})
        =\sum_I\beta_{I,j}\det H_I.
\]
Concavity and homogeneity at \(\Id\) give
\(P_j(K)^{1/k}\le DP_j(\Id):K/k\) for \(K\ge0\).
Use the diagonal matrix with entries \(d_{i,j}^{-1}\) on a
fixed \(I\) and zero elsewhere. This gives \(\beta_{I,j}\le1\).
The defining identities also give
\begin{equation}\label{eq:exactmarginal}
 \sum_{I\ni i}\beta_{I,j}
          \prod_{\ell\in I\setminus\{i\}}d_{\ell,j}=1.
\end{equation}
Since \(0<d_{\ell,j}\le1\), it follows that
\[
 1\le Q_j(\Id)\le\binom mk,\qquad DQ_j(\Id)\ge\Id.
\]
The normalized \(Q_j\) remain in \(\PP_{m,k}\) by
\cref{eq:congruence}. Extract limits of their finitely many
coefficients and of \(d_j\).
Since \(\sum_i d_i^*=k\), \(d_i^*\le1\), and some \(d_i^*=0\),
we have \(k\le a<m\).

The domains \(\mathcal D_j\) contain \(B_1\), their projections
onto \(\R^A\) are uniformly bounded, and they contain segments
of diverging length in each \(B\)-direction.
Take a local Hausdorff limit \(C\) of their closures.
It contains \(B_1\) and the whole subspace \(\R^B\).
A closed convex set containing that subspace is invariant under
translation by it: for \(x\in C\), \(b\in\R^B\), take the limit of
\((1-\varepsilon)x+\varepsilon(b/\varepsilon)\).
Thus the interior of \(C\) is
\(\Omega_*\times\R^B\), with \(B_1^a\subset\Omega_*\subset B_C^a\).
Convexity and the fixed interior ball ensure that each compact
subset of this cylinder lies in \(\mathcal D_j\) for all large \(j\).

Quadratic comparison gives
\(-R^2/2\le w_j\le0\) and \(w_j(0)\le-1/2\).
Local Lipschitz bounds for bounded convex functions therefore
give a locally uniform subsequential limit \(W_\infty\) in the
cylinder. Its restriction to every full \(B\)-parallel line is
bounded and convex, hence constant. Write
\(W_\infty(z)=\varphi(z_A)\).
By local uniform convergence of \(E_{Q_j}\) and
viscosity stability, \(Q_\infty(D^2W_\infty)=1\).

For each \(i\in A\), \cref{eq:exactmarginal} passes to
\begin{equation}\label{eq:activemarginal}
 \sum_{\substack{I\subset A\\i\in I}}\beta_I^*
           \prod_{\ell\in I\setminus\{i\}}d_\ell^*=1.
\end{equation}
Every omitted term contains a vanishing \(B\)-factor.
This proves \(\widehat Q(\Id)>0\) and
\(D\widehat Q(\Id)\ge\Id\).
The restriction assertion follows from \cref{eq:restriction}.
We next verify \cref{eq:activePDE} in the sense of
\cref{def:viscosity}.
Indeed \(\Id_A\oplus0_B\) is a positive semidefinite matrix
where \(Q_\infty>0\), and belongs to \(\Gamma_{Q_\infty}\).
Using it as hyperbolicity direction shows that the intersection
of the full cone with \(\Sym_A\oplus0_B\) is exactly
\(\Gamma_{\widehat Q}\oplus0_B\).
Lifting upper and lower tests of \(\varphi\), independently of
the \(B\)-variables, consequently gives \cref{eq:activePDE}.
This argument includes indefinite lower tests.

For the boundary values, consider
\(\Omega_j^A=\{z_A:(z_A,0)\in\mathcal D_j\}\).
These convex sets are uniformly bounded, contain \(B_1^a\), and
converge to \(\overline\Omega_*\) in Hausdorff distance after taking
closures.
Equation \cref{eq:activemarginal} implies, for all large \(j\),
\[
 \sum_{\substack{I\subset A\\i\in I}}\beta_{I,j}\ge1/2
 \qquad(i\in A).
\]
Every boundary point of this coordinate section lies on the
boundary of \(\mathcal D_j\), where \(W_j=0\).
By \cref{lem:slices}, applied using only \(I\subset A\),
\[
 |W_j(z_A,0)|^k\le C\dist(z_A,\partial\Omega_j^A).
\]
The support functions converge uniformly. For a convex body
with support function \(s\), the boundary distance at an interior
point \(z_A\) equals \(\min_{|e|=1}(s(e)-e\cdot z_A)\).
Passing to the limit, we obtain
\[
 |\varphi(z_A)|^k\le C\dist(z_A,\partial\Omega_*).
\]
Thus \(\varphi\) extends continuously by zero to the boundary.
The bounds in \cref{eq:fullcoverage} follow by coefficient convergence.
\end{proof}

\subsection{Exclusion of nontrivial contacts}

We need two local facts: regularity passes from a smooth limit to
nearby solutions, and a smooth part of a convex solution with zero
boundary values has no flat directions.

\begin{lemma}\label[lemma]{lem:transfer}
Let \(P_j\) be positive multiples of members of \(\PP_{m,k}\),
converging in coefficients to \(P\), with
\(P(\Id)>0\) and \(DP(\Id)>0\).
Suppose convex viscosity solutions \(v_j\) of
\(P_j(D^2v_j)=1\) converge locally uniformly near zero to a
smooth function \(v\) satisfying \(P(D^2v)=1\).
Then, on a fixed neighborhood of zero and for all large \(j\),
the functions \(v_j\) are \(C^{2,\alpha}\) with a uniform bound.
In particular \(D^2v_j(0)\to D^2v(0)\).
\end{lemma}

\begin{proof}
Put \(H=D^2v(0)\ge0\). By \cref{lem:logmono},
\(DP(H)>0\). For all large \(j\), define
\(H_j=P_j(H)^{-1/k}H\). Then \(H_j\to H\) and \(P_j(H_j)=1\).
The functions \(M\mapsto E_{P_j}(H_j+M)\) are
globally degenerate elliptic and vanish at zero. By
\cref{eq:DE} and the implicit function theorem they have
common ellipticity constants and common \(C^2\) bounds in a
fixed neighborhood of zero.
Savin's perturbation theorem \cite{Savin} therefore applies with
a common smallness threshold; only local uniform ellipticity is
needed. Choose \(r\) so that the quadratic Taylor remainder of
\(v\), divided by \(r^2\), is below this threshold.
For large \(j\), the functions
\[
 r^{-2}\big(v_j(rx)-v(0)-rDv(0)\cdot x\big)
             -\tfrac12x^TH_jx
\]
are small in \(C^0(B_1)\) and solve the corresponding shifted
equations. By Savin's estimate, the functions \(v_j\) have a
common \(C^{2,\alpha}\) bound on \(B_{r/2}\). Compactness and the
known \(C^0\) limit identify their Hessian limit as \(D^2v\).
\end{proof}

\begin{lemma}\label[lemma]{lem:smoothpositive}
Let \(P\in\PP_{m,k}\) with \(DP(\Id)>0\).
Let \(v\) be a convex viscosity solution of \(P(D^2v)=1\) on a
bounded convex domain, continuous and zero on its boundary. At every point where \(v\) is smooth
in a neighborhood, \(D^2v\) is positive definite.
\end{lemma}

\begin{proof}
By convexity and \(P(D^2v)=1\), we have \(v<0\) in the domain.
By \cref{lem:jacobi},
\(\{M>0:P(M^{-1})\le1\}\) is convex.
The polynomial \(F=P-1\) satisfies \(F(0)=-1\), and is elliptic
at every matrix \(H\ge0\) with \(P(H)=1\), by
\cref{lem:logmono}. The constant rank theorem
\cite{BG} therefore applies on each smooth neighborhood.
It gives constant rank there and locally fixed kernel directions.

Suppose a smooth point has a kernel direction \(e\).
Subtract its supporting affine function \(\ell\), translate the
point to zero, and write \(g=v-\ell\ge0\).
In a small cylinder, with \(s\) in direction \(e\), we have
\(g(s,z)=\psi(z)\). The function \(\psi\) is smooth and convex.
The supporting plane at \((0,z)\) has no \(s\)-component, so
\[
 g(s,z)\ge\psi(z)
\]
whenever \((s,z)\) belongs to the original domain and \(z\) is
sufficiently small. On a thin cylinder about any compact
subinterval of the \(e\)-axis containing zero, the smooth function
\(\Phi(s,z)=\psi(z)\) solves the same equation.
Its Hessians range over a compact subset of \(\{P=1\}\) on
which \(DP\) is positive definite.

For strong comparison, suppose \(g-\Phi\) is positive somewhere.
Choose a ball \(B_r(y)\subset\{g>\Phi\}\) tangent to \(\{g=\Phi\}\)
at an interior point. Set
\(\zeta=e^{-\alpha|x-y|^2}-e^{-\alpha r^2}\).
On \(r/2<|x-y|<r\), for \(A(x)=DE_P(D^2\Phi(x))\) and large
\(\alpha\),
\[
 e^{-\alpha|x-y|^2}
 \bigl(4\alpha^2 A(x)[x-y,x-y]-2\alpha\tr A(x)\bigr)>0,
\]
because \(A\) has uniform positive lower and upper bounds.
For small \(\varepsilon>0\),
\(E_P(D^2\Phi+\varepsilon D^2\zeta)>0\) and
\(\Phi+\varepsilon\zeta\le g\) on both boundary spheres.
By the supersolution condition, this ordering holds in the annulus.
It also holds outside the ball near the tangency point, where
\(\zeta\le0\). Hence \(\Phi+\varepsilon\zeta\) is a lower test
there with positive operator value, a contradiction.

Hence \(g=\Phi\) on each such cylinder. Exhausting the axial
chord shows that \(v=\ell\) along the whole chord. Its two
endpoints lie on the zero boundary, so \(\ell\), and hence \(v\),
vanishes on that chord. This contradicts \(v<0\).
\end{proof}

\begin{proof}[Proof of \cref{thm:weak}]
For \(k=1\), the equation is linear and uniformly elliptic,
so \(v\) is smooth. Each \(v_{ee}\ge0\) solves the same homogeneous
linear equation. If it vanished at an interior point, it would
vanish everywhere. Then \(v\) would be affine along that direction,
contradicting its zero boundary values and strict interior negativity.

Fix \(k\ge2\) and induct on \(m\ge k\).
If \(m=k\), the class consists only of the determinant.
By the Monge--Amp\`ere localization theorem, a nontrivial contact
set has no interior extreme point \cite{CaffLoc}.
All its extreme points would therefore lie on the zero boundary.
The supporting affine function would vanish on their convex hull,
contradicting \(v<0\) at an interior contact point.
Thus \(v\) is strictly convex. The equivalence of viscosity and
Alexandrov solutions and interior regularity imply
smoothness, and \(\det D^2v=1\) implies \(D^2v>0\).

Let \(m>k\), and assume the theorem in every dimension
\(k\le a<m\). We construct sections near an exposed contact point, following
the localization argument in \cite{CaffLoc}, keeping track of the
transformed coefficients. Suppose a supporting affine function \(\ell\)
has a nontrivial contact set
\[
 K=\{x\in\overline\Omega:v(x)=\ell(x)\}
\]
through an interior point. As above, \(K\) has an interior
extreme point. By Straszewicz's theorem \cite{Rockafellar}, \(K\) has an
interior exposed point: a hyperplane meets \(K\) only at that
point and has all of \(K\) on one side.
Translate this point to zero and put \(g=v-\ell\ge0\).
Choose a unit vector \(\nu\) such that
\begin{equation}\label{eq:exposed}
 K\subset\{\nu\cdot x\ge0\},\qquad
 K\cap\{\nu\cdot x=0\}=\{0\}.
\end{equation}
Choose a nonzero \(q\in K\cap\Omega\), shortening a contact
segment if necessary. Thus \(\nu\cdot q>0\).

A bounded convex body with barycenter \(p\) admits a centered
ellipsoid with
\begin{equation}\label{eq:centeredellipsoid}
 p+LB_1\subset S\subset p+LB_m;
\end{equation}
this follows from the inertia ellipsoid estimate
\cite{KLS}. Consequently, with \(c_m=1/(2m)\),
\[
 x\in\overline S\quad\Longrightarrow\quad
 (1+c_m)p-c_mx\in S.
\]
Choose once and for all
\[
 p_*=tq,\qquad 0<t<c_m/(1+c_m).
\]
If \(K\cap\partial\Omega\ne\varnothing\), its compactness and
\cref{eq:exposed} give
\(a_0=\min_{K\cap\partial\Omega}\nu\cdot x>0\).
Choose \(t\) still smaller so that
\begin{equation}\label{eq:centerchoice}
 (1+c_m)\nu\cdot p_*<c_ma_0.
\end{equation}
We keep \(p_*\) fixed as the section height tends to zero.

For each \(h>0\), minimize over \(p\in\R^m\) the convex function
\[
 \mathcal F_h(p)=\int_\Omega
          (h-g(x)+p\cdot(x-p_*))_+\dd.
\]
The function \(\mathcal F_h\) is finite and continuous.
Integrating on a ball about \(p_*\), we obtain
\(\mathcal F_h(p)\ge c|p|-C\), so it has a minimizer \(p_h\).
The level \(g-p_h\cdot(x-p_*)=h\) has zero \(m\)-dimensional
measure, since it is above its value zero at \(p_*\).
Differentiating \(\mathcal F_h\) at \(p_h\), we find that the
bounded convex section
\[
 S_h=\{x\in\Omega:g(x)<h+p_h\cdot(x-p_*)\}
\]
has barycenter \(p_*\).
At first it may touch \(\partial\Omega\).
By \cref{eq:centeredellipsoid}, for
\(x\in\overline S_h\),
\[
 0\le g(x)\le(1+c_m^{-1})h.
\]
Indeed apply its affine upper bound to
\((1+c_m)p_*-c_mx\), where \(g\ge0\).
Every cluster point of these sections as \(h\downarrow0\)
therefore lies in \(K\).
In fact the sections are uniformly compactly contained in
\(\Omega\) for small \(h\).
Otherwise \(x_h\in\overline S_h\) would approach a point
\(x\in K\cap\partial\Omega\), while their balanced images would
approach a point of \(K\) with negative \(\nu\)-coordinate,
by \cref{eq:centerchoice}. If \(K\Subset\Omega\), the conclusion
follows because every limit point of the sections belongs to \(K\).

For all small \(h\), \(p_h\cdot q<0\).
If not, \(q\in S_h\) along a sequence, so the fixed point
\[
 (1+c_m)p_*-c_mq=(t-c_m(1-t))q
\]
would be a limiting contact point with negative \(\nu\)-coordinate.
It follows that \(0\in S_h\) and its section depth there is
\begin{equation}\label{eq:contactdepth}
 h-p_h\cdot p_*\ge h.
\end{equation}
Let \(b_h=-\min_{\overline S_h}\nu\cdot x\).
Then \(b_h>0\) and \(b_h\to0\).

Normalize \(S_h\) as in \cref{eq:centeredellipsoid}, centered
at \(p_*\), and put
\[
 s_h=P(L_h^{-T}L_h^{-1})^{-1/k},\qquad
 w_h(y)=\frac{g(p_*+L_hy)-h-p_h\cdot L_hy}{s_h}.
\]
The domain of \(w_h\) lies between \(B_1\) and \(B_m\), and
\(w_h\) has zero boundary values. Define
\[
 P_h(H)=s_h^kP(L_h^{-T}HL_h^{-1}).
\]
Then \(P_h(\Id)=1\) and \(P_h(D^2w_h)=1\).
By \cref{eq:congruence,lem:logmono}, we have
\(P_h\in\PP_{m,k}\) and \(DP_h(\Id)>0\).
Quadratic comparison gives
\[
 -m^2/2\le w_h\le0,\qquad
 h/s_h=-w_h(0)\in[1/2,m^2/2].
\]
At \(y_h=-L_h^{-1}p_*\), we have
\(-w_h(y_h)\ge1/2\) by \cref{eq:contactdepth}.
The distance of \(y_h\) to the supporting plane
\(\nu\cdot(p_*+L_hy)=-b_h\) is \(b_h/|L_h^T\nu|\).
Since \(|y_h|\le m\),
\[
 |L_h^T\nu|\ge(\nu\cdot p_*)/m>0.
\]
Thus this distance tends to zero.
Set \(e_h=L_h^T\nu/|L_h^T\nu|\).
By \cref{lem:slices},
\[
 e_h^TDP_h(\Id)e_h\to0.
\]
In particular the smallest eigenvalue of \(DP_h(\Id)\) tends
to zero.

Apply \cref{lem:newton}, absorbing its orthogonal coordinates
in \(L_h\). The number \(a\) of remaining variables satisfies \(k\le a<m\).
Its limit \(\varphi\) is convex, continuous and zero on the
boundary of a bounded convex domain, and solves
\(\widehat Q(D^2\varphi)=1\) with \(D\widehat Q(\Id)>0\).
Let \(Z=\widehat Q(\Id)\).
The operator \(\widehat Q/Z\) and the function
\(Z^{1/k}\varphi\) satisfy the normalized induction hypothesis.
Hence \(\varphi\) is smooth and \(D^2\varphi>0\).
The full limit \(W_\infty(z)=\varphi(z_A)\) is smooth as well.
Its Hessian at zero is \(H_*=D^2\varphi(0)\oplus0_B\);
\cref{eq:fullcoverage} and \cref{lem:logmono} give
\[
 DQ_\infty(H_*)\ge
 \frac{1}{\binom mk\,\tr H_*}\Id>0.
\]
By \cref{lem:transfer} and successive elliptic estimates, \(W_h\) is
smooth near zero for every sufficiently small \(h\) in the chosen
sequence.
Fix one such \(h\). Both affine changes are invertible, so
the original \(v\) is smooth near the fixed point \(p_*\).
It is affine along the contact segment through \(p_*\),
so its Hessian there has a kernel.
This contradicts \cref{lem:smoothpositive}.

Every supporting affine function at an interior point therefore
touches \(v\) only at that point.
\Cref{prop:strictcontact} gives smoothness and positive definite
Hessian throughout \(\Omega\), completing the induction.
\end{proof}

\section{Entire solutions}\label{sec:entire}

\begin{proof}[Proof of \cref{thm:main}]
Replacing \(u\) by \(c_0^{-1/k}u\), we may assume \(c_0=1\).
For \(k=1\), the equation is \(\Delta u=1\). Convexity bounds
every entry of \(D^2u\); these entries are harmonic and hence
constant. We therefore assume \(k\ge2\).

At each point, \(H=D^2u\ge0\) has rank at least \(k\).
The linearization of \(\sigma_k\) is positive definite there.
Indeed its diagonal entries in an eigenbasis are
\(\sigma_{k-1}(\lambda\mid i)>0\).
Continuity of \(D^2u\), local ellipticity and concavity of
\(\sigma_k^{1/k}\) give local Evans--Krylov estimates and
then smoothness. Applying the constant-rank theorem of
Bian and Guan \cite{BG} as in \cref{lem:smoothpositive}, now without a
boundary argument, shows that the Hessian has constant rank
\(r\) and a fixed kernel. The local kernels agree on overlaps,
so the kernel is fixed on \(\R^n\). After an orthogonal change
of variables and subtraction of an affine function,
\[
 u(x',x'')=v(x'),\qquad x'\in\R^r,\qquad
 D^2v>0,\qquad \sigma_k(D^2v)=1.
\]
Here \(r\ge k\). If \(r=k\), the J\"orgens--Calabi--Pogorelov
theorem \cite{Jorgens,Calabi,Pogorelov} applies.
It remains to consider the positive definite case \(k<n\),
renaming \(r\) as \(n\). We fix the original point zero and
subtract its tangent plane, so that
\begin{equation}\label{eq:fixedpoint}
 u(0)=0,\qquad Du(0)=0,\qquad H_0=D^2u(0)>0.
\end{equation}

For every \(h>0\), there is \(p_h\in\R^n\) such that
\[
 K_h=\{x:u(x)<p_h\cdot x+h\}
\]
is bounded and has centroid zero.
To construct the centered section, minimize the extended-valued
convex function
\[
 \mathcal F_h(p)=\int_{\R^n}(h-u(x)+p\cdot x)_+\,dx.
\]
Since \(D^2u(0)>0\) and \(u\) is convex, \(u(x)\ge c|x|\)
for large \(|x|\), so \(\mathcal F_h(0)<\infty\).
Integrating on a fixed ball, we obtain
\(\mathcal F_h(p)\ge c_h|p|-C_h\).
Lower semicontinuity follows from Fatou's lemma, hence a finite
minimizer exists.
If its positive set were unbounded, that open convex set,
which contains zero, would contain a ray from zero.
Writing \(f(x)=u(x)-p\cdot x\), convexity and \(f(0)=0\)
would give \(f\le0\) along the ray. Consequently
\(\{f<h/2\}\) would contain a ball and a ray, and have
infinite volume. This contradicts finiteness of the integral.
Thus the positive set is bounded. Beyond a larger sphere, \(f\)
has a uniform linear lower bound by convexity. The positive sets
for nearby slopes therefore lie in a common ball.
Since \(\{f=h\}\) has measure zero, we may differentiate at the
minimizer to obtain \(\int_{K_h}x\,dx=0\).

Choose a centered ellipsoid by \cite{KLS}:
\[
 L_hB_1\subset K_h\subset L_hB_R,\qquad R=n.
\]
We may write \(L_h=O_h\diag(\ell_{1,h},\ldots,\ell_{n,h})\),
where \(O_h\) is orthogonal and \(\ell_{i,h}>0\).
Suppress the index \(h\) and put \(B_i=h/\ell_i^2\).
Set \(\theta=\sigma_k(B)^{1/k}\), \(b_i=B_i/\theta\), and
\(D_b=\diag(b_i)\). On \(\Omega_h=L_h^{-1}K_h\), define
\[
 w_h(y)=\frac\theta h
       \bigl(u(L_hy)-p_h\cdot L_hy-h\bigr).
\]
The corresponding operator is
\[
 P_h(M)=\sigma_k(D_b^{1/2}MD_b^{1/2}).
\]
Then \(P_h(\Id)=\sigma_k(b)=1\),
\(P_h(D^2w_h)=1\), and \(w_h=0\) on \(\partial\Omega_h\).
Since \(w_h(0)=-\theta\), by comparison with quadratic barriers
on \(B_1\) and \(B_R\),
\begin{equation}\label{eq:theta}
 \frac12\le\theta\le\frac{R^2}{2},\qquad
 -\frac{R^2}{2}\le w_h\le0.
\end{equation}
Introduce the Newton weights
\[
 T_i=\sigma_{k-1}(b\mid i),\qquad
 d_i=b_iT_i,\qquad D=\diag(d_i).
\]
They satisfy \(0<d_i\le1\), \(\sum_i d_i=k\), and
\(DP_h(\Id)=D\).

We now prove that these weights have a uniform positive lower
bound as \(h\to\infty\). More precisely, for all sufficiently
large \(h\),
\begin{equation}\label{eq:noncollapse}
 d_i\ge\eta>0\quad \text{for }1\le i\le n,
\end{equation}
where \(\eta\) may depend on \(u\).
Otherwise there is a sequence \(h\to\infty\) along which
the smallest weight tends to zero. Apply \cref{lem:newton}
to \(w_h\) and \(P_h\), after taking a subsequence and
permuting coordinates. For
\(W_h(z)=w_h(D^{1/2}z)\), the limit is
\(W_\infty(z)=\varphi(z_A)\) on
\(\Omega_*\times\R^B\), with \(B\ne\varnothing\).
The function \(\varphi\) solves
\(\widehat Q(D^2\varphi)=1\) on \(\Omega_*\), with zero boundary
values. By \cref{thm:weak}, after normalizing \(\widehat Q\),
it is smooth and has positive definite Hessian.
The full limiting operator is elliptic at
\(D^2W_\infty(0)\), by \cref{eq:fullcoverage,lem:logmono}.
By \cref{lem:transfer},
\[
 D^2W_h(0)\longrightarrow D^2\varphi(0)\oplus0_B.
\]
On the other hand, both normalizations preserve the original
point zero. Differentiating their definitions, we obtain
\begin{equation}\label{eq:fixedidentity}
 D^2W_h(0)=\diag(\sqrt{T_i})\,
          O_h^TH_0O_h\,\diag(\sqrt{T_i}).
\end{equation}
For \(i\in B\),
\(\sigma_k(b\mid i)=1-d_i\to1\). Since \(k\le n-1\),
the Newton--Maclaurin inequality in the remaining \(n-1\)
variables implies
\[
 T_i\ge
 \binom{n-1}{k-1}
 \left(\frac{1-d_i}{\binom{n-1}{k}}\right)^{(k-1)/k}.
\]
By \cref{eq:fixedidentity},
\[
 (W_h)_{ii}(0)\ge\lambda_{\min}(H_0)T_i\ge c>0.
\]
This contradicts the vanishing \(B\)-block of the limiting Hessian
and proves \cref{eq:noncollapse}.

It remains to convert \cref{eq:noncollapse} into a global
Hessian bound. By \cref{thm:classical},
\begin{equation}\label{eq:whbounds}
 c\Id\le D^2w_h\le C\Id\quad\hbox{on }B_{1/2},
\end{equation}
uniformly for all large \(h\).
At zero we also have
\[
 D^2w_h(0)=D_b^{-1/2}O_h^TH_0O_hD_b^{-1/2}.
\]
By comparing diagonal entries and using \cref{eq:whbounds},
we bound every \(b_i\) above and below. Together with
\cref{eq:theta}, we obtain \(c\sqrt h\le\ell_i\le C\sqrt h\).
Hence \(L_hB_{1/2}\) contains \(B_{c\sqrt h}\).
Transform \cref{eq:whbounds} back to the original variables
and let \(h\to\infty\). We obtain
\begin{equation}\label{eq:global}
 c\Id\le D^2u(x)\le C\Id\qquad(x\in\R^n).
\end{equation}

On the compact set of matrices in \cref{eq:global},
\(F=\sigma_k^{1/k}\) is concave and uniformly elliptic.
Taking the infimum of its tangent affine functions at those
matrices extends it to a concave uniformly elliptic operator
on all of \(\Sym_n\), agreeing with \(F\) on that set.
The functions \(u_R(x)=R^{-2}u(Rx)\) have uniformly bounded
\(C^{1,1}(B_1)\) norm, by \cref{eq:fixedpoint,eq:global}.
By the Evans--Krylov estimate \cite{CS}, for fixed \(x,y\)
and all sufficiently large \(R\),
\[
 |D^2u(x)-D^2u(y)|\le C\left(\frac{|x-y|}{R}\right)^\alpha.
\]
Letting \(R\to\infty\) proves that \(D^2u\) is constant.
Restoring the discarded affine function, the flat directions,
and the factor \(c_0^{1/k}\) proves the theorem.
\end{proof}

\section*{AI Disclosure}
The approach taken in this paper grew out of the authors' earlier
work on Monge--Amp\`ere equations and their use of Caffarelli's
section and localization methods \cite{CaffLoc}. The strategy combines
section normalization, estimates in the original and dual
variables, and dimension reduction to control degeneration.
The authors used chatGPT 5.6 Sol and 6 Astra to do the matrix
calculation for the inverse Hessian Jacobi inequality in
\cref{lem:jacobi}, the logarithmic integral estimate
\eqref{eq:logL1}, and the double logarithm and contact set
argument in the proof of \cref{thm:classical}. AI was also used to work out the details of the dual
approximation in \cref{prop:strictcontact} and the rescaling
and boundary estimates in \cref{lem:slices,lem:newton}.
It further helped derive and check the Hessian identity
\eqref{eq:fixedidentity} and the inequalities leading to
\eqref{eq:noncollapse}.
AI further assisted with checking proof dependencies,
locating references, and revising the exposition.
The authors are responsible for the mathematical content,
its verification, and the final manuscript.


\begin{thebibliography}{99}
\raggedright

\bibitem{AD}
G.~Acosta and R.~G.~Dur\'an,
An optimal Poincar\'e inequality in \(L^1\) for convex domains,
\emph{Proc. Amer. Math. Soc.} \textbf{132} (2004), 195--202.


\bibitem{BCGJ}
J.~Bao, J.~Chen, B.~Guan, and M.~Ji,
Liouville property and regularity of a Hessian quotient equation,
\emph{Amer. J. Math.} \textbf{125} (2003), 301--316.

\bibitem{BGLS}
H.~H.~Bauschke, O.~G\"uler, A.~S.~Lewis, and H.~S.~Sendov,
Hyperbolic polynomials and convex analysis,
\emph{Canad. J. Math.} \textbf{53} (2001), 470--488.


\bibitem{BG}
B.~Bian and P.~Guan,
A structural condition for microscopic convexity principle,
\emph{Discrete Contin. Dyn. Syst.} \textbf{28} (2010), 789--807.


\bibitem{BBL}
J.~Borcea, P.~Br\"and\'en, and T.~M.~Liggett,
Negative dependence and the geometry of polynomials,
\emph{J. Amer. Math. Soc.} \textbf{22} (2009), 521--567.

\bibitem{CaffLoc}
L.~A.~Caffarelli,
A localization property of viscosity solutions to the Monge--Amp\`ere
equation and their strict convexity,
\emph{Ann. of Math. (2)} \textbf{131} (1990), 129--134.

\bibitem{CNS}
L.~Caffarelli, L.~Nirenberg, and J.~Spruck,
The Dirichlet problem for nonlinear second order elliptic equations,
III: Functions of the eigenvalues of the Hessian,
\emph{Acta Math.} \textbf{155} (1985), 261--301.

\bibitem{CS}
L.~Caffarelli and L.~Silvestre,
On the Evans--Krylov theorem,
\emph{Proc. Amer. Math. Soc.} \textbf{138} (2010), 263--265.

\bibitem{Calabi}
E.~Calabi,
Improper affine hyperspheres of convex type and a generalization
of a theorem by K.~J\"orgens,
\emph{Michigan Math. J.} \textbf{5} (1958), 105--126.

\bibitem{CY}
S.-Y.~A.~Chang and Y.~Yuan,
A Liouville problem for the Sigma-2 equation,
\emph{Discrete Contin. Dyn. Syst.} \textbf{28} (2010), 659--664.

\bibitem{CJTZ}
R.~Chen, H.~Jian, X.~Tu, and X.~Zhou,
Regularity for convex viscosity solutions of \(\sigma_2\) equation,
preprint (2026).
\href{https://arxiv.org/abs/2605.30823v2}{arXiv:2605.30823v2}.


\bibitem{ChengYau}
S.-Y.~Cheng and S.-T.~Yau,
Complete affine hypersurfaces. Part I. The completeness of affine metrics,
\emph{Comm. Pure Appl. Math.} \textbf{39} (1986), 839--866.

\bibitem{CW}
K.-S.~Chou and X.-J.~Wang,
A variational theory of the Hessian equation,
\emph{Comm. Pure Appl. Math.} \textbf{54} (2001), 1029--1064.

\bibitem{CD}
J.~Chu and S.~Dinew,
Liouville theorem for a class of Hessian equations,
preprint (2023).
\href{https://arxiv.org/abs/2306.13825v1}{arXiv:2306.13825v1}.

\bibitem{CIL}
M.~G.~Crandall, H.~Ishii, and P.-L.~Lions,
User's guide to viscosity solutions of second order partial
differential equations,
\emph{Bull. Amer. Math. Soc.} \textbf{27} (1992), 1--67.

\bibitem{Du}
S.-Z.~Du,
Necessary and sufficient conditions to Bernstein theorem of a Hessian equation,
\emph{Trans. Amer. Math. Soc.} \textbf{375} (2022), no.~7, 4873--4892.

\bibitem{GS}
P.~Guan and M.~Sroka,
A special concavity property for positive Hessian quotient operators,
\emph{Discrete Contin. Dyn. Syst.} \textbf{54} (2026), 50--60.


\bibitem{Hu}
X.~Hu,
Sharp Hausdorff bounds for the interior singular set of convex
\(k\)-Hessian solutions,
preprint (2026).
\href{https://arxiv.org/abs/2607.28988v1}{arXiv:2607.28988v1}.


\bibitem{Jorgens}
K.~J\"orgens,
\"Uber die L\"osungen der Differentialgleichung \(rt-s^2=1\),
\emph{Math. Ann.} \textbf{127} (1954), 130--134.

\bibitem{KLS}
R.~Kannan, L.~Lov\'asz, and M.~Simonovits,
Isoperimetric problems for convex bodies and a localization lemma,
\emph{Discrete Comput. Geom.} \textbf{13} (1995), 541--559.

\bibitem{LRW}
M.~Li, C.~Ren, and Z.~Wang,
An interior estimate for convex solutions and a rigidity theorem,
\emph{J. Funct. Anal.} \textbf{270} (2016), 2691--2714.

\bibitem{Mooney25}
C.~Mooney,
Remarks on the quadratic Hessian equation,
\emph{J. Math. Study}, to appear.
\href{https://arxiv.org/abs/2505.14586v1}{arXiv:2505.14586v1}.

\bibitem{Moreau}
J.-J.~Moreau,
Proximit\'e et dualit\'e dans un espace hilbertien,
\emph{Bull. Soc. Math. France} \textbf{93} (1965), 273--299.

\bibitem{Pogorelov}
A.~V.~Pogorelov,
On the improper convex affine hyperspheres,
\emph{Geom. Dedicata} \textbf{1} (1972), 33--46.

\bibitem{Rockafellar}
R.~T.~Rockafellar,
\emph{Convex Analysis}, Princeton Mathematical Series, vol.~28,
Princeton University Press, Princeton, NJ, 1970.

\bibitem{Savin}
O.~Savin,
Small perturbation solutions for elliptic equations,
\emph{Comm. Partial Differential Equations} \textbf{32} (2007), 557--578.

\bibitem{SavinObstacle}
O.~Savin,
The obstacle problem for Monge Ampere equation,
\emph{Calc. Var. Partial Differential Equations} \textbf{22} (2005), 303--320.

\bibitem{SY}
R.~Shankar and Y.~Yuan,
Rigidity for general semiconvex entire solutions to the sigma-2 equation,
\emph{Duke Math. J.} \textbf{171} (2022), 3201--3214.

\bibitem{Trudinger}
N.~S.~Trudinger,
On the Dirichlet problem for Hessian equations,
\emph{Acta Math.} \textbf{175} (1995), 151--164.

\bibitem{Tu}
Q.~Tu,
Pogorelov type estimates for \((n-1)\)-Hessian equations and related
rigidity theorems,
preprint (2024).
\href{https://arxiv.org/abs/2405.02939v2}{arXiv:2405.02939v2}.

\bibitem{BW}
B.~Wang,
A Liouville type rigidity result for entire solutions of the \(k\)-Hessian equation,
preprint (2026).
\href{https://arxiv.org/abs/2608.29392v2}{arXiv:2608.29392v2}.

\bibitem{WB}
C.~Wang and J.~Bao,
Liouville property and existence of entire solutions of Hessian equations,
\emph{Nonlinear Anal.} \textbf{223} (2022), 113020.

\bibitem{Warren}
M.~Warren,
Nonpolynomial entire solutions to \(\sigma_k\) equations,
\emph{Comm. Partial Differential Equations} \textbf{41} (2016), 848--853.

\bibitem{Zhang}
R.~Zhang,
\(C^2\) estimates for \(k\)-Hessian equations and a rigidity theorem,
\emph{Adv. Math.} \textbf{480} (2025), 110488.

\end{thebibliography}
\end{document}